%%%%%plain TeX file of a bastard.tex, an article by D. Zeilberger
 
%begin macros

\baselineskip=14pt
\parskip=10pt
\def\Tilde{\char126\relax}
\def\halmos{\hbox{\vrule height0.15cm width0.01cm\vbox{\hrule height
 0.01cm width0.2cm \vskip0.15cm \hrule height 0.01cm width0.2cm}\vrule
 height0.15cm width 0.01cm}}
\font\eightrm=cmr8  
\font\eighttt=cmtt8
\magnification=\magstephalf

\parindent=0pt
\overfullrule=0in
%end macros
\bf
\centerline{How Much Should a 19th-Century French Bastard Inherit}
\medskip
\it
\centerline{ Doron ZEILBERGER \footnote{$^1$}
{\baselineskip=9pt
\eightrm  \raggedright
Department of Mathematics, Temple University,
Philadelphia, PA 19122, USA. 
{\baselineskip=9pt
\eighttt zeilberg@math.temple.edu ; 
\break http://www.math.temple.edu/\Tilde zeilberg }.
Supported in part by the NSF. Feb. 6, 1997.
}}
\medskip
\qquad\qquad\qquad\qquad\qquad\qquad\qquad
\qquad\qquad
{\it Dedicated to Gerry Ladas on his 60th Birthday.}
\medskip
\rm
{\bf Abstract:}
Catalan's formula, for the portion of the inheritance that
a legitimate child of 
a 19th-century deceased French gentleman should receive, is
given a new proof (using Difference Operators), and generalized.
Another, more computationally efficient, formula is also derived.
 
{\bf Catalan's Formula}
 
Chrystal's classical text[Ch] (II, p. 416) contains the following
charming exercise, due to Catalan[Ca].
 
{\it By French law an illegitimate child receives one-third
of the portion of  the inheritance that he would have received had he
been legitimate. If there are $l$ legitimate and $n$ illegitimate
children, show that the portion of inheritance $1$ due to a legitimate
child is}
$$
{{1} \over {l}} -
{{n} \over {3l(l+1)}} +{{n(n-1)} \over {3^2l(l+1)(l+2)}} -
\dots
\pm{{n(n-1) \dots 2.1} \over {3^n l(l+1)\dots (l+n)}} \quad .
$$
 
Let's call the portion due to a legitimate child $a(l,n)$.
If the condition would have been that
each illegitimate child receives one-third of what a legitimate
child is to receive, this would have yielded the simple algebraic equation
$la(l,n)+n(1/3)a(l,n)=1$ and 
hence $a(l,n)=3/(3l+n)$. But since the bastard is to receive
one-third of {\it what he would have received had he been legitimate},
we have a {\it difference} (or recurrence) equation:
$$
la(l,n)+(n/3)a(l+1,n-1)=1 \quad .
\eqno(*)
$$
 
Introducing the shift operators
$Lf(l,n):=f(l+1,n)$ and $Nf(l,n)=f(l,n+1)$, this can be
rewritten:
$$
(1+{{n} \over {3l}}L N^{-1})a(l,n)={{1} \over {l}} \quad,
$$
yielding:
$$
a(l,n)=(1+{{n} \over {3l}}L N^{-1})^{-1} {{1} \over {l}}=
\sum_{r=0}^{\infty} (-{{n} \over {3l}}L N^{-1})^r{{1} \over {l}}
$$
$$
=\sum_{r=0}^{\infty}
(-1)^r 
({{n} \over {3l}}{{n-1} \over {3(l+1)}} \dots{{n-r+1} \over {3(l+r-1)}}
L^r N^{-r}){{1} \over {l}}=
\sum_{r=0}^{n}
{(-1)^r {n(n-1) \dots (n-r+1)} \over {3^r l(l+1) \dots (l+r)}} \quad. \halmos
$$
 
\eject
{\bf Another Formula}
 
While Catalan's formula is elegant, it is not any better than
using the recurrence $(*)$ directly, by starting with
$a(l+n,0)=1/(l+n)$ and successively computing $a(l+n-1,1)$,
$a(l+n-2,2), \dots, a(l,n)$, requiring $O(n)$ operations
and $O(1)$ storage. The recurrence $(*)$ enables us to
quickly find $a(n+1,l-1)$ or $a(n-1,l+1)$, once we know
$a(n,l)$, in case the status of one of the children changes.
However, neither the recurrence $(*)$, nor Catalan's formula, are
efficient if we have already computed $a(l,n)$, and all of a sudden
another illegitimate child, that was previously unknown to the family,
shows up at the lawyer's office. 
By using {\tt zeillim} in the package {\tt EKHAD} that accompanies
[PWZ], we find the recurrence
$$
a(l,n+1)={{1} \over {l+n+1}}+{{2(n+1)} \over {3(l+n+1)}}a(l,n) \quad,
$$
that implies the alternative formula:
$$
a(l,n)={{n! (2/3)^n} \over {(n+l)!}} \sum_{b=0}^{n}
{{(l+b-1)!} \over {b!}} (3/2)^b \quad.
$$
 
{\bf The Multi-Mistress Case}
 
A typical well-to-do 19th-century French gentleman must have 
had {\it plusieurs amantes}.
It is very unlikely that he liked them equally well.
Let's say that he had $m$ of them, and he wished to leave to 
every one of the $n_i$ children of Mistress $i$, $x_i$ of what he would
have left them had they been legitimate ($i=1, \dots, m$).
Suppose that he also had $l$ legitimate children. What portion of
the inheritance should each of these $l$ legitimate children receive?
Let's call this quantity $a(l; n_1, \dots , n_m)$. Then,
$$
\left ( 1+\sum_{i=1}^{m}
{{n_i}x_i \over {l}}L N_i^{-1} \right )
a(l;n_1, \dots , n_m)={{1} \over {l}} \quad,
\quad and \quad hence
$$
$$
a(l;n_1, \dots , n_m)=
$$
$$
\left ( 1+\sum_{i=1}^{m}
{{n_i}x_i \over {l}}L N_i^{-1} \right )^{-1} 
{{1} \over {l}}=
\sum_{r_1, \dots , r_m \geq 0}^{\infty} 
(-1)^{r_1+ \dots + r_m} {{(r_1+ \dots + r_m)!} \over {r_1! \cdots r_m!}}
\left ( \prod_{j=1}^{m} ({{x_j n_j} \over {l}}L N_j^{-1})^{r_j} \right ) 
{{1} \over {l}}
$$
$$
=\sum_{r_1, \dots , r_m \geq 0}^{\infty} 
(-1)^{r_1+ \dots + r_m} {{(r_1+ \dots + r_m)!} \over {r_1! \cdots r_m!}}
\left ( \prod_{j=1}^{m} x_j^{r_j} \right ) 
\left ( \prod_{j=1}^{m} (n_j N_j^{-1})^{r_j} \right ) 
((1/l)L)^{r_1+ \dots + r_m} {{1} \over {l}}
$$
$$
=
\sum_{r_1, \dots , r_m \geq 0}^{\infty} 
(-1)^{r_1+ \dots + r_m} {{(r_1+ \dots + r_m)!} \over {r_1! \cdots r_m!}}
{{1} \over {(l)(l+1) \cdots (l+r_1 + \dots + r_j)}}
\prod_{j=1}^{m} x_j^{r_j} 
\prod_{j=1}^{m} {{n_j!} \over {(n_j-r_j)!}}
 \quad.
$$
Hence
$$
a(l; n_1, \dots , n_m)=
\sum_{r_1=0}^{n_1} \dots \sum_{r_m=0}^{n_m}
{{\prod_{j=1}^{m} (-x_j)^{r_j} {{n_j} \choose {r_j}}}
\over 
{ (1+r_1+ \dots + r_m) {{l+ r_1 + \dots + r_m} \choose {l-1}} } } \quad 
. \halmos
$$
\eject
{\bf REFERENCES}
 
[Ca] E.C. Catalan, {\it Nouv. Ann.}, ser. II, t.2.
 
[Ch] G. Chrystal, {\it ``ALGEBRA''}, reprinted by Chelsea, New York, NY,
1964.
 
[PWZ] M. Petkovsek, H.S. Wilf, and D. Zeilberger, {\it ``A=B''},
A.K.Peters, Wellesley, MA, 1996. The accompanying Maple package EKHAD may be
downloaded from the author's website.
 
\bye